\documentclass[12pt,a4paper]{amsart}
\usepackage{tikz}
\usepackage{pdflscape}
\usepackage{amsmath}
\usepackage{amsfonts}
\usepackage{amssymb}
\usepackage{mathtools}
\usepackage[all]{xy}
\usepackage{color}
\usetikzlibrary{arrows}
\usepackage{float}
\usepackage[shortlabels]{enumitem}
\usepackage{ifpdf}
\ifpdf
  \usepackage[colorlinks,
            linkcolor=magenta,       
            anchorcolor=magenta,     
            citecolor=magenta,       
            final,
            hyperindex]{hyperref}
\else
  \usepackage[colorlinks,
            linkcolor=magenta,       
            anchorcolor=magenta,     
            citecolor=magenta,       
            final,
            hyperindex,driverfallback=hypertex]{hyperref}
\fi
\newtheorem{theorem}{Theorem}[section]
\newtheorem{lemma}[theorem]{Lemma}
\newtheorem{corollary}[theorem]{Corollary}

\newtheorem{proposition}[theorem]{Proposition}
\theoremstyle{definition}
\newtheorem{defn}[theorem]{Definition}

\newtheorem{example}[theorem]{Example}

\makeatletter \@addtoreset{equation}{section} \makeatother
\allowdisplaybreaks

\newcommand{\delete}[1]{}

\newcommand{\marg}{\sim_{\omega}}
\newcommand{\cg}{\equiv_{\omega}}
\newcommand{\HH}{\mathcal H}

\newcommand{\CC}{\mathcal C}

\newcommand{\omone}{\omega_1}
\newcommand{\omtwo}{\omega_2}

\begin{document}
\title[Marginal Congruences and Elements]{Marginal Congruences and Elements for Two Commutator Words in Completely Regular Semigroups}

\author{Lin Shen  and  Shoufeng Wang}
\address{School of Mathematics, Yunnan Normal University, Kunming, Yunnan 650500, China}
\email{wsf1004@163.com}


\begin{abstract}Motivated by Problem 13 posed by Ara\'ujo, Ara\'ujo
and Kinyon \cite{Araujo2025} concerning marginal
elements, subsemigroups and marginal equivalences associated
with commutator words in completely regular semigroups,
in this paper we investigate the two commutator words $
\omega_1(x,y)=x^{-1}y^{-1}xy$ and $
\omega_2(x,y)=(yx)^{-1}xy$ on completely regular semigroups.
We prove that their marginal equivalences  $\sim_{\omega_r}$ ($r=1,2$)
coincide and are congruences. More precisely, for a completely regular semigroup $S$,  we show that  $a\sim_{\omega_r}b$ ($r=1,2$) if and only if
$a^0=b^0=p$ and $b^{-1}a\in Z(pSp)$,  where $Z(pSp)$ is the center of the monoid $pSp$, and so  both words have the same  set of marginal elements:
\[
M_{\omega}(S)
=\bigcup_{p\in E(S)}
\bigl(Z(pSp)\cap H_p\bigr),
\] where $E(S)$ is the set of idempotents in $S$.
We also characterize when this set is a subsemigroup
and establish alternative descriptions of the common
marginal congruence and marginal elements.
\end{abstract}
\makeatletter
\@namedef{subjclassname@2020}{\textup{2020} Mathematics Subject Classification}
\makeatother

\maketitle



\vspace{-0.8cm}

\section{Introduction}\label{sec:introduction}
The notions of verbal and marginal subgroups associated
with group words were introduced by Hall \cite{Hall1940}.
Completely regular semigroups, which are unions of groups,
form an important class of semigroups whose structure has
been extensively studied; see
\cite{Clifford1961,Petrich1999}.
The study of commutator words in this setting is more
subtle than in groups, since the inverse of a product
need not equal the product of the inverses in reverse order.

Recent investigations of commutators in completely simple
semigroups \cite{Radovic2025} and of abelianness and
centrality in inverse semigroups \cite{Kinyon2026}
provide further motivation for studying commutator
structures in semigroups.
In particular, Ara\'ujo, Ara\'ujo and Kinyon
\cite{Araujo2025} investigated marginal subsemigroups
and commutators in inverse semigroups.
In \cite[Problem~13]{Araujo2025}, they raised the problem
of studying marginal elements, marginal subsemigroups
and marginal congruences associated with commutator
words in completely regular semigroups.

We consider the commutator words $\omega_1(x,y)=x^{-1}y^{-1}xy$ and
$\omega_2(x,y)=(yx)^{-1}xy$ on completely regular semigroups.
For cryptogroups, Ara\'ujo, Ara\'ujo and Kinyon \cite{Araujo2025} showed
that the marginal equivalence for each word is a congruence.
For all completely regular semigroups, we prove that the two marginal
equivalences coincide and form a congruence. We characterize this
congruence via centers of local monoids, determine the common marginal
elements, and give a criterion for their set to be a subsemigroup.
We also give alternative characterizations of the congruence and
the marginal elements.

\section{Main results}
\label{sec:CR}

\label{sec:regular-defs}

In this section, $S$ is always a completely regular semigroup. Every $a\in S$
belongs to its maximal subgroup $H_{a^0}$, where $a^0=aa^{-1}=a^{-1}a$ and $a^{-1}$ is called the {\em local inverse} of $a$. A \emph{completely regular semigroup word} is a
term in multiplication and local inverse.

\begin{defn}\label{def:marg} Let $S$ be a completely regular semigroup.
Assume that $\omega(x_1,\ldots,x_n)$ is a completely regular semigroup word.
For $a,b\in S$, write $a\marg b$ if and only if the following two identities hold
for every $i\in\{1,\ldots,n\}$ and every $(x_1,\ldots,x_n)\in S^n$:
\begin{align}
 \omega(x_1,\ldots,ax_i,\ldots,x_n)
 &=\omega(x_1,\ldots,bx_i,\ldots,x_n), \label{eq:mleft}\\
 \omega(x_1,\ldots,x_i a,\ldots,x_n)
 &=\omega(x_1,\ldots,x_i b,\ldots,x_n). \label{eq:mright}
\end{align}
The relation $\marg$ is always an equivalence relation. Let $\cg$ be the smallest congruence containing $\marg$.
Then
\begin{equation}\label{eq:congruence-criterion}
 \marg=\cg
 \quad\mbox{if and only if}\quad
 \marg\text{ is  a congruence}.
\end{equation}
Denote
\begin{equation}\label{eq:canonical}
M_\omega(S)=\{a\in S:a\sim_\omega a^0\}
\end{equation}
and call the elements of $M_\omega(S)$ {\em marginal elements of $S$ for the word $\omega(x_1,\ldots,x_n)$.}
\end{defn}

In what follows, we consider two commutator words:
\begin{equation}\label{eq:two-words}
 \omone(x,y)=x^{-1}y^{-1}xy,
 \qquad
 \omtwo(x,y)=(yx)^{-1}xy.
\end{equation}
Unlike the group and inverse-semigroup settings,
these words may differ because
$(yx)^{-1}$ need not be $x^{-1}y^{-1}$.
For either word, marginal equivalence requires all four
identities below to hold for all $x,y\in S$:
\begin{align}
 \omega(ax,y)&=\omega(bx,y), \tag{M1}\\
 \omega(xa,y)&=\omega(xb,y), \tag{M2}\\
 \omega(x,ay)&=\omega(x,by), \tag{M3}\\
 \omega(x,ya)&=\omega(x,yb). \tag{M4}
\end{align}
We first consider the case of completely simple semigroups. Recall that  a completely simple semigroup is exactly a Rees matrix semigroup over a group.
Let $T=M[G;I,\Lambda;P]$ be a completely simple semigroup  with the multiplication
\begin{equation}\label{eq:reesmult}
 (i,g,\lambda)(j,h,\mu)
 =(i,gp_{\lambda j}h,\mu).
\end{equation}
Observe that the inverse of $(i,g,\lambda)$ within its maximal subgroup is
\[
 (i,g,\lambda)^{-1}
 =\bigl(i,p_{\lambda i}^{-1}g^{-1}
             p_{\lambda i}^{-1},\lambda\bigr).
\]
and the identity of $H_{i\lambda}$ is
$e_{i\lambda}=(i,p_{\lambda i}^{-1},\lambda)$.

\begin{proposition}\label{thm:rees}
Let $T=M[G;I,\Lambda;P]$ be a completely simple  semigroup and $Z(G)$ be the center of $G$. For $r=1,2$,
\begin{align}
 (i,g,\lambda)\sim_{\omega_r}(j,h,\mu)
 \quad\Longleftrightarrow\quad
 i=j,\ \lambda=\mu,\ gh^{-1}\in Z(G).
 \label{eq:reesrelation}
\end{align}
In particular, $\sim_{\omega_1}=\sim_{\omega_2}$ is a
congruence, with quotient
$
 T/\!\sim_{\omega_r}
 \cong M[G/Z(G);I,\Lambda;\overline P],
$
where $\overline P=(p_{\lambda i}Z(G))_{\lambda\in\Lambda,\,i\in I}$.
Moreover,
\begin{equation}\label{eq:reesmarg}
 M_{\omega_r}(T)
  =\{(i,g,\lambda):gp_{\lambda i}\in Z(G)\}
  =\bigcup_{e\in E(T)}Z(H_e).
\end{equation}
So $M_{\omega_r}(T)$ is a subsemigroup if and only if
$p_{\lambda i}^{-1}p_{\lambda j}
 p_{\mu j}^{-1}p_{\mu i}\in Z(G)
$  for all
$i,j\in I$ and $\lambda,\mu\in\Lambda$.
\end{proposition}
\begin{proof}
Let $
u=(i,g,\lambda), v=(j,h,\mu),
$
and put
\[
A=p_{\lambda i},\quad B=p_{\mu j},\quad
C=p_{\lambda j},\quad D=p_{\mu i}.
\]
Then
$
u^{-1}=(i,A^{-1}g^{-1}A^{-1},\lambda),
$
and a direct calculation gives
\begin{align}
\omega_1(u,v)
&=(i,A^{-1}g^{-1}A^{-1}CB^{-1}
h^{-1}B^{-1}DgCh,\mu), \label{eq:omega1-Rees}\\
\omega_2(u,v)
&=(j,C^{-1}g^{-1}D^{-1}h^{-1}
C^{-1}AgCh,\mu). \label{eq:omega2-Rees}
\end{align}

We first prove the necessity of the stated conditions.
Suppose that $u\sim_{\omega_r}v$.
For $\omega_1$, the row index of $\omega_1(x,y)$
is determined by the first argument, whereas its
column index is determined by the second argument.
Consequently, the marginal identities
\[
\omega_1(ux,y)=\omega_1(vx,y)
\quad\text{and}\quad
\omega_1(x,yu)=\omega_1(x,yv)
\]
imply, respectively, that
$i=j$ and $\lambda=\mu.$
For $\omega_2$, both the row and column indices of
$\omega_2(x,y)$ are determined by the second
argument. Thus the identities
\[
\omega_2(x,uy)=\omega_2(x,vy)
\quad\text{and}\quad
\omega_2(x,yu)=\omega_2(x,yv)
\]
give the same conclusions.
Hence, in either case, $u$ and $v$ belong to the
same maximal subgroup $H_{i\lambda}$.
Restricting all substitutions to $H_{i\lambda}$,
both $\omega_1$ and $\omega_2$ reduce to the
ordinary group commutator. In particular, taking
$x=v^{-1}$ in the left marginal identity yields
\[
[u v^{-1},y]=[vv^{-1},y]=1_{H_{i\lambda}}
\qquad (y\in H_{i\lambda}).
\]
It follows that
$
uv^{-1}\in Z(H_{i\lambda}).
$
Under the natural group isomorphism
\[
H_{i\lambda}\longrightarrow G,\qquad
(i,k,\lambda)\longmapsto k p_{\lambda i},
\]
this condition is equivalent to the fact that
$
gh^{-1}\in Z(G).
$

We now prove the converse. Suppose that
$
i=j, \lambda=\mu,
gh^{-1}\in Z(G).
$
Put $z=gh^{-1}\in Z(G)$, so that $g=zh$.
We claim that, for either $r=1,2$, replacing
the group coordinate of either argument of
$\omega_r$ by its product with a central element
does not change the value of $\omega_r$.
Indeed, this follows directly from
\eqref{eq:omega1-Rees} and \eqref{eq:omega2-Rees}:
in each expression, the group coordinate of
each argument occurs once with exponent $1$
and once with exponent $-1$, and the central
factor can be moved across all sandwich entries
and cancelled. More explicitly, for arbitrary
$x=(k,t,\nu)\in T$, we have
\[
ux=(i,gp_{\lambda k}t,\nu),\qquad
vx=(i,hp_{\lambda k}t,\nu),
\]
whose group coordinates differ by the central
factor $z$. Similarly,
\[
xu=(k,tp_{\nu i}g,\lambda),\qquad
xv=(k,tp_{\nu i}h,\lambda),
\]
and their group coordinates also differ by $z$.
The same holds for the pairs $(uy,vy)$ and $(yu,yv)$ for every $y\in T$.

Therefore, by the coordinate invariance just
established, we obtain
\begin{align*}
\omega_r(ux,y)=\omega_r(vx,y),\,\,
\omega_r(xu,y)=\omega_r(xv,y),\\
\omega_r(x,uy)=\omega_r(x,vy),\,\,
\omega_r(x,yu)=\omega_r(x,yv)
\end{align*}
for all $x,y\in T$ and $r=1,2$.
Thus $u\sim_{\omega_r}v$, proving the desired
equivalence.

The quotient map
\[
\varphi:T\longrightarrow M[G/Z(G);I,\Lambda;\overline{P}],
\qquad
(i,g,\lambda)\longmapsto(i,gZ(G),\lambda),
\]
is a surjective homomorphism, and its kernel
is precisely $\sim_{\omega_r}$. Hence
$\sim_{\omega_r}$ is a congruence and
\[
T/\mathord{\sim_{\omega_r}}
\cong M[G/Z(G);I,\Lambda;\overline{P}].
\]

Finally, the identity of $H_{i\lambda}$ is
$
e_{i\lambda}=(i,p_{\lambda i}^{-1},\lambda).
$
Applying the characterization above to
$u\sim_{\omega_r}e_{i\lambda}$, we obtain
\[
(i,g,\lambda)\in M_{\omega_r}(T)
\iff gp_{\lambda i}\in Z(G).
\]
Consequently,
\[
M_{\omega_r}(T)
=\{(i,g,\lambda)\in T:
gp_{\lambda i}\in Z(G)\}
=\bigcup_{e\in E(T)}Z(H_e).
\]

To determine when this set is a subsemigroup,
write two arbitrary elements of $M_{\omega_r}(T)$
as
$
a=(i,zp_{\lambda i}^{-1},\lambda)$ and $
b=(j,wp_{\mu j}^{-1},\mu),
$
where $z,w\in Z(G)$. Their product is
$
ab=(i,zp_{\lambda i}^{-1}p_{\lambda j}
wp_{\mu j}^{-1},\mu).
$
Thus $ab\in M_{\omega_r}(T)$ if and only if
$
zp_{\lambda i}^{-1}p_{\lambda j}
wp_{\mu j}^{-1}p_{\mu i}\in Z(G),
$
which, since $z,w\in Z(G)$, is equivalent to
$
p_{\lambda i}^{-1}p_{\lambda j}
p_{\mu j}^{-1}p_{\mu i}\in Z(G).
$
Therefore $M_{\omega_r}(T)$ is a subsemigroup
if and only if this condition holds for all
$i,j\in I$ and $\lambda,\mu\in\Lambda$.
\end{proof}

The following example gives a completely simple semigroup $T$ in which  $M_{\omega_r}(T)$  is not a subsemigroup.
\begin{example}
\label{ex:24}
Take $G=S_3$, $I=\Lambda=\{1,2\}$ and
\[
 P=\begin{pmatrix}1&1\\[2pt]1&t\end{pmatrix},
 \qquad t=(12)\in S_3.
\]
Since $Z(S_3)=\{1\}$, we have $M_{\omega_r}(T)=E(T)$ by
\eqref{eq:reesmarg}.
Obviously,
$ e=(1,1,1), f=(2,t,2)\in E(T)$.  But
$ ef=(1,t,2)\notin E(T).$
\end{example}

In the following statements, we shall prove  that $\sim_{\omega_1}=\sim_{\omega_2}$ and it is a
congruence for general completely regular semigroups.  Let $S$ be a completely regular semigroup and $p\in E(S)$. Then $pSp$ is
a monoid with identity $p$. Let $Z(pSp)$ be  the ordinary center of $pSp$.
Recall that $H_p$ is the maximal subgroup of $S$ containing $p$.
We observe that,
if $q\in E(S), t\in S$, $x\in H_q$ and $t x=x$ (respectively $x t=x$), then
$tq=q$ (respectively $qt=q$).
Define a relation $\mathrel{\CC}$ as follows: For all $a,b\in S$,
\begin{equation}\label{eq:C}
 a\mathrel{\CC}b
 \quad\Longleftrightarrow\quad
 a^0=b^0\quad\text{and}\quad b^{-1}a\in Z(a^0Sa^0).
\end{equation}

\begin{lemma}\label{lem:corner}
Suppose $p,q\in E(S)$ and $z\in Z(pSp)\cap H_p$.
\begin{enumerate}
\item If $pq=q$, then $zq\in Z(qSq)\cap H_q$,
      with $(zq)^{-1}=z^{-1}q$.
\item If $qp=q$, then $qz\in Z(qSq)\cap H_q$,
      with $(qz)^{-1}=qz^{-1}$.
\end{enumerate}
\end{lemma}
\begin{proof}
Assume $pq=q$.  The element $qp$ is an idempotent
of $pSp$.  Hence centrality of $z$ gives
$
 zqp=qpz=qz.
$
Multiplying on the right by $q$ gives
$
 zq=qzq.
$
Since $z^{-1}$ is also a central unit of $pSp$, we likewise have
$z^{-1}q=qz^{-1}q$.
Thus $q$ is a two-sided identity for both $zq$ and $z^{-1}q$, and
$
 (zq)(z^{-1}q)=z(qz^{-1}q)=zz^{-1}q=pq=q.
$
Interchanging $z$ and $z^{-1}$ gives $(z^{-1}q)(zq)=q$.
Thus $zq\in H_q$, with the asserted inverse.
If $t\in qSq$, then $pt=t$ and $tq=t$,
while $tp\in pSp$.  Thus
$
 ztp=tpz=tz.
$
Multiplying on the right by $q$ yields
$zt=tzq$.  Consequently
$
 (zq)t=zt=tzq=t(zq),
$
and $zq$ is central in $qSq$.
The second assertion is the dual statement,
obtained in the opposite semigroup.
\end{proof}

\begin{lemma}\label{prop:C}
The relation $\CC$ is a congruence.
\end{lemma}
\begin{proof}
For each idempotent $p$, the group
$K_p=Z(pSp)\cap H_p$ is a central subgroup of $H_p$.
The restriction of $\CC$ to $H_p$ is the coset relation
modulo $K_p$, and different maximal subgroups are
disjoint.  Hence $\CC$ is an equivalence relation.
Let $a\mathrel{\CC}b$, with $a^0=b^0=p$.
Put $z=b^{-1}a\in K_p$.  Then $a=bz=zb.$
For any $t\in S$, let $r=(bt)^0$.
Since $p(bt)=bt$, we have $pr=r$.
By Lemma~\ref{lem:corner}, $c=zr\in Z(rSr)\cap H_r$.  Hence
$
 at=bzt=z(bt)=(zr)(bt)=c(bt)=(bt)c.
$
In particular, $(at)^0=(bt)^0=r$ and
$ (bt)^{-1}(at)=c\in Z(rSr).$
Thus $at\mathrel{\CC}bt$.
Similarly, let $s=(tb)^0$.  Since $(tb)p=tb$,
one has $sp=s$.  The dual part of
Lemma~\ref{lem:corner} gives
$d=sz\in Z(sSs)\cap H_s$, and
$ ta=(tb)z=(tb)(sz)=(tb)d.$
It follows that $ta\mathrel{\CC}tb$.
Therefore $\CC$ is a semigroup congruence.
Lastly, $a=bz$ inside $H_p$ gives
$a^{-1}=z^{-1}b^{-1}$.  Since $z$ is central
in $H_p$ and in $pSp$, we have
$(b^{-1})^{-1}a^{-1}=ba^{-1}=z^{-1}\in K_p$.
Thus $a^{-1}\mathrel{\CC}b^{-1}$.
\end{proof}

\begin{lemma}\label{lem:invariance}
For $i=1,2$, if $u\mathrel{\CC}u'$ and
$v\mathrel{\CC}v'$, then
$
 \omega_i(u,v)=\omega_i(u',v').
$
\end{lemma}
\begin{proof}
It suffices to replace one argument at a time.
Write $e=u^0$ and $f=v^0$.
First take $u'=uc$ with
$c\in Z(eSe)\cap H_e$.  Since $uc\in H_e$
and $c$ commutes with $u$, we have
\[
\begin{aligned}
 \omega_1(uc,v)
 =c^{-1}(u^{-1}v^{-1}u)c\,v
 =u^{-1}v^{-1}uv=\omega_1(u,v),
\end{aligned}
\]
as $u^{-1}v^{-1}u\in eSe$.
Next let $v'=vd$, with
$d\in Z(fSf)\cap H_f$.  Likewise,
\[
\begin{aligned}
 \omega_1(u,vd)
 =u^{-1}d^{-1}(v^{-1}uv)d
 =u^{-1}v^{-1}uv=\omega_1(u,v),
\end{aligned}
\]
since $v^{-1}uv\in fSf$.

For $\omega_2$, again take $u'=uc$.
Put $w=vu$ and $q=w^0$. Since $we=w$, we have
$q e=q$.  Applying the \emph{second} part of
Lemma~\ref{lem:corner} gives
$
 k=qc\in Z(qSq)\cap H_q$ and $ k^{-1}=qc^{-1}.
$
Because $u$ commutes with $c$ and $wq=w$, we obtain that
$
 v(uc)=wc=w(qc)=wk$ and $
 (uc)v=c(uv).
$
The elements $wk,w$ both belong to $H_q$, so
$(wk)^{-1}=k^{-1}w^{-1}$.  Moreover,
$w^{-1}q=w^{-1}$ and $k$ commutes with
$w^{-1}\in qSq$.  Hence, \emph{without assuming that
$uv$ belongs to $H_q$}, we have
\[
\begin{aligned}
 \omega_2(uc,v)
 =k^{-1}w^{-1}c(uv)
 =k^{-1}w^{-1}qc(uv)
 =k^{-1}w^{-1}k(uv)
 =w^{-1}(uv)=\omega_2(u,v).
\end{aligned}
\]
The insertion of $q$ uses $w^{-1}q=w^{-1}$.
Since $k$ commutes with $w^{-1}$, the cancellation follows from
$k^{-1}w^{-1}k=k^{-1}k\,w^{-1}=q\,w^{-1}=w^{-1}$, where
$q$ is the identity of $H_q$.

Finally, take $v'=vd$, and again set
$r=(vu)^0$.  Since $f(vu)=vu$, we have
$fr=r$.  The \emph{first} part of
Lemma~\ref{lem:corner} shows that
$dr\in Z(rSr)\cap H_r$ and
$(dr)^{-1}=d^{-1}r$.
As $v d=d v$, we have
\[
\begin{aligned}
 (vdu)^{-1}
 =(d(vu))^{-1}
 =((dr)(vu))^{-1}
=(vu)^{-1}(dr)^{-1}
 =(dr)^{-1}(vu)^{-1}
 =d^{-1}(vu)^{-1}.
\end{aligned}
\]
Therefore
\[
 \omega_2(u,vd)
 =d^{-1}\underbrace{(vu)^{-1}uv}_{W}d.
\]
But $W\in fSf$: indeed
$f(vu)^{-1}=(vu)^{-1}$ and $(uv)f=uv$.
Centrality of $d$ in $fSf$ now yields
$d^{-1}Wd=W$, as required.
\end{proof}

\begin{corollary}\label{cor:subset}
For each $i\in\{1,2\}$,
$\CC\subseteq\sim_{\omega_i}$.
\end{corollary}
\begin{proof}
If $a\mathrel{\CC}b$, then, by
Lemma~\ref{prop:C}, all four pairs
$(ax,bx)$, $(xa,xb)$, $(ay,by)$ and
$(ya,yb)$ lie in $\CC$.
Apply Lemma~\ref{lem:invariance} to each
replacement in the four marginal identities (M1)--(M4).
\end{proof}

\begin{lemma}\label{lem:sameH}
If $a\sim_{\omega_i}b$, then $a\mathrel{\HH}b$ and,
writing $p=a^0=b^0$, then
$
 z=ab^{-1}=b^{-1}a \in Z(H_p).
$
\end{lemma}
\begin{proof}
A completely regular semigroup is a semilattice of
completely simple $\mathcal D$-classes.  If $a,b$
lie in components indexed by $\alpha,\beta$, then
both words project to the meet of their variable indices.
Testing the first marginal equation with $x,y$ in
the $\alpha$-component, then the $\beta$-component,
gives $\alpha=\alpha\wedge\beta=\beta$.
Thus $a,b$ lie in a single completely simple component.
Restricting the testing variables to that component,
Proposition~\ref{thm:rees} implies $a\mathrel{\HH}b$.
Within $H_p$ both words reduce to the ordinary group
commutator.  Taking $x=b^{-1}$ and $y\in H_p$ in
the first marginal identity gives
$
 [ab^{-1},y]=[bb^{-1},y]=p.
$
Hence $ab^{-1}\in Z(H_p)$; commuting with $b$ shows
$ab^{-1}=b^{-1}a$.
\end{proof}

\begin{lemma}
\label{lem:detect}
Let $T$ be a completely regular monoid with
identity $p$, and let $z\in H_p$ be a unit.
For either $\omega_1$ or $\omega_2$, if
$z\sim_{\omega_i}p$ in $T$, then $z\in Z(T)$.
\end{lemma}
\begin{proof}
For both words and every $x\in T$, we have
$ \omega_i(x,p)=x^0.$
Applying the first two marginal equations
with second argument $p$, we get
$
 (zx)^0=(xz)^0=x^0
$ for all $x\in T$.
In particular, if $q\in E(T)$, then
$zq,qz\in H_q$.  Since $zq,qz\in H_q$, both have $q$ as a left
and right identity.  Hence
$$ zq=q(zq)=qzq=(qz)q=qz.$$
Now take an arbitrary $x\in T$, and put
$q=x^0$.  We have
$zx,xz\in H_q$ and
$zq=qz\in H_q$.

If $i=2$, the marginal equation
$\omega_2(z,x)=\omega_2(p,x)$ gives
$
 (xz)^{-1}(zx)=q
$
inside the group $H_q$, and hence $zx=xz$.

If $i=1$, the analogous equation yields
$
 z^{-1}x^{-1}zx=q.
$
Because $z$ commutes with $q$, this is the
ordinary group commutator identity
$
 (zq)^{-1}x^{-1}(zq)x=q
  \text{in }H_q.
$
Thus $zq$ commutes with $x$; consequently
$zx=(zq)x=x(zq)=xz$.
In both cases $z$ is central in $T$.
\end{proof}

\begin{lemma}\label{prop:reverse}
For each $i\in\{1,2\}$,
$\sim_{\omega_i}\subseteq\CC$.
\end{lemma}
\begin{proof}
Assume $a\sim_{\omega_i}b$.
By Lemma~\ref{lem:sameH}, we have
$a,b\in H_p$ and
$
 z=ab^{-1}=b^{-1}a\in Z(H_p).
$
We claim that $z\sim_{\omega_i}p$ in $S$.
Indeed, substituting $x=b^{-1}t$ into the
first marginal equation for $a,b$ gives
$
 \omega_i(zt,y)=\omega_i(pt,y).
$
Substituting $x=tb^{-1}$ into the second
gives
$
 \omega_i(tz,y)=\omega_i(tp,y).
$  The corresponding substitutions
$y=b^{-1}t$ and $y=tb^{-1}$ into the
third and fourth equations give
$
\omega_i(x,zt)=\omega_i(x,pt)$ and
$\omega_i(x,tz)=\omega_i(x,tp).
$
These are exactly the four marginal
identities for $z$ and $p$.
Now restrict all testing variables to
$T=pSp$.  This is a completely regular
submonoid with identity $p$: if
$x\in pSp$, then $px=x=xp$ and hence
$px^{-1}=x^{-1}=x^{-1}p$.
The element $z$ is a unit of $T$, and
the restricted marginal identities say
$z\sim_{\omega_i}p$ in $T$.
By Lemma~\ref{lem:detect},
$z\in Z(T)=Z(pSp)$.
This is precisely $a\mathrel{\CC}b$.
\end{proof}
\begin{theorem}\label{thm:main}
For every completely regular semigroup $S$, we have
$$
\sim_{\omega_1}
 =\sim_{\omega_2}
 =\equiv_{\omega_1}
 =\equiv_{\omega_2}
 =\CC.
$$
\end{theorem}
\begin{proof}
Lemma~\ref{prop:reverse} and
Corollary~\ref{cor:subset} give
$
 \sim_{\omega_1}=\CC=\sim_{\omega_2}.
$
By Lemma~\ref{prop:C}, $\CC$ is
already a congruence, including
compatibility with local inversion.
Consequently its generated congruence
is itself, so the result holds.
\end{proof}

\begin{corollary}
\label{cor:common-margin}
For every completely regular semigroup $S$,
\[
 M_{\omega_1}(S)=M_{\omega_2}(S)
 =\bigcup_{p\in E(S)}(Z(pSp)\cap H_p).
 \]
\end{corollary}
\begin{proof}
An element $a\in H_p$ belongs to $ M_{\omega_r}(S)$ if and only if  $a\sim_{\omega_r}p$ for $r=1,2$.
By Theorem~\ref{thm:main} and (\ref{eq:C}), this is equivalent to
$a =p^{-1}a\in Z(pSp)$.
\end{proof}

The following example shows that $Z(pSp)$ can not be changed to $Z(H_p)$  in Corollary \ref{cor:common-margin}.
\begin{example}
\label{ex:B-two-level}
Let $C_2$ be written additively as $\{0,1\}$.
Take the two-element semilattice of groups with
upper group $G=C_2\times C_2$, lower group
$L=S_3\times C_2$, and connecting homomorphism
\[
 \phi:G\longrightarrow L,
 \qquad \phi(s,t)=\bigl((12)^s,t\bigr).
\]
Then we have a Clifford semigroup $S=G\bigcup L$. Let $p=(0,0)$ be the
identity in the upper group; it is the global identity,
so $pSp=S$ and $H_p=G$.
Every element of $G$ lies in the group center $Z(H_p)$,
since $G$ is abelian. But for $a=(1,0)\in G$, its
image $\phi(a)=((12),0)$ does not commute with the
lower-group element $((123),0)$.
Accordingly $a\notin Z(pSp)$, and
\[
 Z(H_p)=G\quad\text{whereas}\quad
 Z(pSp)\cap H_p=\{(0,0),(0,1)\}.
\]
Thus replacing $Z(pSp)$ by $Z(H_p)$ would introduce
false marginal elements, even in a Clifford semigroup.
\end{example}

The following result gives a sufficient and necessary condition under which  $M_{\omega_1}(S)$ and $M_{\omega_2}(S)$ are subsemigroups of $S$.
\begin{corollary}
\label{cor:B-idempotent-closure}Let $S$ be  a completely regular semigroup and denote $$K=\bigcup\limits_{p\in E(S)}(Z(pSp)\cap H_p).$$ Then $K$ is a subsemigroup of $S$ if and only if
$ef\in K$ for all $e,f\in E(S)$.
\end{corollary}
\begin{proof}
Let $\rho:S\to S/\mathcal C$ be the quotient map, where
$\mathcal C$ is a congruence by
Lemma~\ref{prop:C}. For an arbitrary $a\in H_p$, we have
 $a^2,a\in  H_p$  and
\[
 \rho(a)^2=\rho(a)
 \ \Longleftrightarrow\ a^2\mathrel{\mathcal C}a
 \ \Longleftrightarrow\ a=a^{-1}a^2\in Z(pSp).
\]
Consequently $K=\rho^{-1}\bigl(E(S/\mathcal C)\bigr).$
 If $K$ is closed, then it contains every product $ef$ for $e,f\in E(S)$
because $E(S)\subseteq K$.   Conversely, assume that the given condition holds and take $a,b\in K$.
By Theorem~\ref{thm:main},
$a\mathrel{\mathcal C}a^0$ and
$b\mathrel{\mathcal C}b^0$, and so
$ab\mathrel{\mathcal C}a^0b^0$. By hypothesis, $a^0b^0\in K$, and so
  $\rho(a^0b^0)$ is idempotent in $S/\mathcal C$ by the fact that $K=\rho^{-1}\bigl(E(S/\mathcal C)\bigr).$   Therefore
$\rho(ab)$ is idempotent as well, and once more by the fact that $K=\rho^{-1}\bigl(E(S/\mathcal C)\bigr)$, we have
  $ab\in K$.
\end{proof}

We next give an alternative characterization of  $\CC$. To this aim, for $i=1,2$, we  define
\begin{equation}\label{eq:B-eta-definition}
 a\mathrel{\eta_i}b
 \quad\Longleftrightarrow\quad
 \omega_i(a,y)=\omega_i(b,y)
 \text{ and }\omega_i(y,a)=\omega_i(y,b)
 \quad\text{for all }y\in S.
\end{equation}
We first prove a useful lemma.

\begin{lemma}
\label{lem:B-pointwise-local}
Let $T$ be a completely regular monoid with identity
$1_T$, and let $a,b$ be units of $T$. If
$a\mathrel{\eta_i}b$ for either commutator,
then $b^{-1}a\in Z(T)$.  For $\omega_2$ alone,
the weaker single-direction hypothesis
\begin{equation}\label{eq:B-one-sided-point}
 \omega_2(y,a)=\omega_2(y,1_T)
 \qquad(y\in T)
\end{equation}
already implies $a\in Z(T)$.
\end{lemma}
\begin{proof}
Write $c=b^{-1}a$. We use the standard description
of a completely regular semigroup as a semilattice
of completely simple components.  Inside any fixed
such component, left or right multiplication by a
fixed element of that component is injective on any
specified maximal subgroup: in Rees coordinates it
multiplies the group coordinate by fixed invertible
sandwich and group factors.  All applications of this
cancellation below take place in one component.

First suppose $i=1$.  For $e\in E(T)$ put
$ f=a^{-1}ea$ and $ g=b^{-1}eb.$
Both are idempotents, since conjugation by a unit is
an automorphism.  The point-value comparisons
$\omega_1(a,e)=\omega_1(b,e)$ and
$\omega_1(e,a)=\omega_1(e,b)$ become
\begin{equation}\label{eq:B-point-A-idempotent}
 fe=ge,\qquad ef=eg.
\end{equation}
A unit lies in the top semilattice component of the
monoid; therefore $e,f,g$ all belong to the same
completely simple component.  The first equality
in \eqref{eq:B-point-A-idempotent} forces $f$ and $g$
to have the same Rees row; the second forces the same
Rees column.  Thus $f\mathrel{\mathcal H}g$.
Two idempotents in one maximal subgroup are equal,
so $f=g$.  As $e$ ranges over $E(T)$, the idempotent
$b^{-1}eb$ does too.  The identity
$a^{-1}ea=b^{-1}eb$ therefore says
$c^{-1}hc=h$ for every $h\in E(T)$; equivalently,
$c$ commutes with every idempotent.

Now compare $\omega_1(y,a)=\omega_1(y,b)$ for arbitrary
$y\in T$. Multiplying both sides on the left by $y$
and writing
$ v=b^{-1}yb and q=y^0$, this gives
\begin{equation}\label{eq:B-corner-cancel-A}
 q(c^{-1}vc)=qv.
\end{equation}
Since $c$ centralizes the idempotents,
$c^{-1}vc$ belongs to the same maximal subgroup
as $v$; the conjugating unit preserves $v^0$.
Moreover $q$ and $v$ lie in the same completely
simple component.  Left multiplication by $q$
is injective on the maximal subgroup of $v$.
Cancelling in \eqref{eq:B-corner-cancel-A} gives
$c^{-1}vc=v$.  Since conjugation by the unit $b$
is a bijection of $T$, the element $v=b^{-1}yb$
ranges over all of $T$.  Thus $c\in Z(T)$.

Next suppose $i=2$. Restricting the two point-value
identities to the unit group $U(T)$ turns
$\omega_2(x,y)$ into the ordinary group commutator.
The group calculation yields
$c\in Z(U(T))$, in particular $a=bc=cb$.
For $e\in E(T)$ the equality
$\omega_2(e,a)=\omega_2(e,b)$ reads
\begin{equation}\label{eq:B-point-B-idempotent}
 (ae)^{-1}ea=(be)^{-1}eb.
\end{equation}
In a completely simple component, the output row
of either side comes from $ae$ or $be$, respectively,
and its output column comes from $ea$ or $eb$.
Therefore \eqref{eq:B-point-B-idempotent} implies
$ ae\mathrel{\mathcal R}be $ and $ ea\mathrel{\mathcal L}eb.
$
Multiplication on the left by a unit preserves
$\mathcal L$-classes, while multiplication on the right
by a unit preserves $\mathcal R$-classes.  Hence
$ae\mathrel{\mathcal L}be$,
$ea\mathrel{\mathcal R}eb$  and
$
 ae\mathrel{\mathcal H}be$, $ ea\mathrel{\mathcal H}eb.
$
Multiplication by a unit on either side preserves $\mathcal H$.
Since $c=b^{-1}a=ab^{-1}$, multiplying the first relation on the left by
$b^{-1}$ and the second on the right by $b^{-1}$ gives
$ce\mathrel{\mathcal H}e$ and $ec\mathrel{\mathcal H}e$, respectively.
Because $e$ is the identity of its maximal subgroup,
the first relation gives $ce=e(ce)=ece$, and the
second gives $ec=(ec)e=ece$.  Hence $ce=ec$ for
every idempotent $e$.

With $y\in T$, put $u=yb$ and $v=by$.
These two elements lie in the same completely simple
component.  Since $c$ commutes with every idempotent,
$uc$ belongs to the maximal subgroup of $u$ and
$ (uc)^{-1}=c^{-1}u^{-1}.$
This inverse identity is justified inside the local
subgroup: with $q=u^0$, the elements $qc$ and $qc^{-1}$
are inverse units of $qTq$, and
$uc=u(qc)$.  The point-value equality
$\omega_2(a,y)=\omega_2(b,y)$ consequently becomes
\begin{equation}\label{eq:B-corner-cancel-B}
 (c^{-1}u^{-1}c)v=u^{-1}v.
\end{equation}
The two left factors belong to the same maximal
subgroup $H_{u^0}$, and multiplication by $v$ is
injective on this subgroup, since $u$ and $v$ lie
in the same completely simple component. Therefore
$c^{-1}u^{-1}c=u^{-1}$.  As $u=yb$ varies over
all of $T$, this proves $c\in Z(T)$.

Finally, suppose only \eqref{eq:B-one-sided-point}.
For an idempotent $e$, it states
$
 (ae)^{-1}ea=e.
$
Comparing Rees rows and columns gives
$ae\mathrel{\mathcal R}e$ and
$ea\mathrel{\mathcal L}e$.  The unit translations
also give $ae\mathrel{\mathcal L}e$ and
$ea\mathrel{\mathcal R}e$.  Hence both $ae$ and $ea$
belong to $H_e$, so
$ae=eae=ea$.  Thus $a$ commutes with all idempotents.
The valid local inverse formula now yields
\[
 \omega_2(y,a)=(ay)^{-1}ya=y^{-1}a^{-1}ya,
 \qquad \omega_2(y,1_T)=y^0.
\]
The assumed equality is therefore
$y^{-1}a^{-1}ya=y^0$.  Multiply on the left by $y$, we have
$
 y^0 a^{-1}ya=y.
$
As $a$ commutes with $y^0$ and $y^0y=y$, the left
side equals $a^{-1}ya$.  Hence $a^{-1}ya=y$
for every $y\in T$, or $a\in Z(T)$.
\end{proof}

\begin{theorem}
\label{thm:B-pointwise}
For every completely regular semigroup $S$, we have
\[
\eta_1=\eta_2=\mathcal C
   =\sim_{\omega_1}=\sim_{\omega_2}
   =\equiv_{\omega_1}=\equiv_{\omega_2}.
\]
\end{theorem}
\begin{proof}
Suppose $a\mathrel{\eta_i}b$ and let $\alpha,\beta$
be the semilattice indices of the completely simple
components containing $a,b$, respectively.
Each $\omega_i$ has, in the semilattice quotient,
the meet of the indices of its two arguments.
Comparing $\omega_i(a,y)$ and $\omega_i(b,y)$ with
$y$ first chosen in the $\alpha$-component and then
in the $\beta$-component gives
$
 \alpha=\alpha\wedge\beta=\beta.
$
Thus $a,b$ lie in the same completely simple component.
For $\omega_1$, equality of $\omega_1(a,y)$ and
$\omega_1(b,y)$ determines the row of $a,b$, while
equality of $\omega_1(y,a)$ and $\omega_1(y,b)$
determines the column. For $\omega_2$, the single
family $\omega_2(y,a)=\omega_2(y,b)$ determines
both row and column.  This follows directly from the
Rees-coordinate formulas in Section~\ref{sec:CR}.
In either case $a\mathrel{\mathcal H}b$, so
$a^0=b^0=p$.
The local monoid $T=pSp$ is completely regular, as
proved in Lemma~\ref{prop:reverse}, and has
identity $p$. Both $a,b$ are units of $T$.
Restrict the equalities defining $\eta_i$ to $y\in T$.
Lemma~\ref{lem:B-pointwise-local} then gives
$b^{-1}a\in Z(T)=Z(pSp)$.  By
\eqref{eq:C}, $a\mathrel{\mathcal C}b$.
This proves $\eta_i\subseteq\mathcal C$.

Conversely, if $a\mathrel{\mathcal C}b$,
Lemma~\ref{lem:invariance}, with the other
argument paired with itself, gives
$\omega_i(a,y)=\omega_i(b,y)$ and
$\omega_i(y,a)=\omega_i(y,b)$ for every $y\in S$.
Hence $\mathcal C\subseteq\eta_i$.
The remaining equalities are exactly those of
Theorem~\ref{thm:main}.
\end{proof}

To give new characterizations of marginal elements, we need the following notion.
\label{subsec:B-existential}
\begin{defn}[c.f. the statements before Theorem 1 in \cite{Araujo2025}]
\label{def:existential}Let $S$ be a  completely regular semigroup and   $\omega(x_1,\ldots,x_m)$ be a  completely regular semigroup word.
For $i\in\{1,\ldots,m\}$, an element $a\in S$ is called
\emph{existentially partially marginal in position $i$} if there
are $e_i,f_i\in E(S)$ such that the following two identities hold
for every $(x_1,\ldots,x_m)\in S^m$:
\[
\begin{aligned}
\omega(\ldots,ax_i,\ldots)
 &=\omega(\ldots,e_i x_i,\ldots),\\
\omega(\ldots,x_i a,\ldots)
 &=\omega(\ldots,x_i f_i,\ldots).
\end{aligned}
\]
Write $\widetilde M_\omega^i(S)$ for this set and
$\widetilde M_\omega(S)=\bigcap_i\widetilde M_\omega^i(S)$.
\end{defn}

\begin{theorem}
\label{thm:B-existential}
Let $S$ be any completely regular semigroup and put
\[
 K=\bigcup_{p\in E(S)}\bigl(Z(pSp)\cap H_p\bigr).
\]
Then
\begin{equation}\label{eq:B-existential-equals}
\widetilde M_{\omega_1}(S)
 =\widetilde M_{\omega_2}(S)
 =M_{\omega_1}(S)=M_{\omega_2}(S)=K.
\end{equation}
\end{theorem}
\begin{proof}
Let $a\in\widetilde M_{\omega_i}(S)$ and set
$p=a^0$.  We first note that every idempotent witness
for any one of its defining insertion identities must
lie in the same completely simple component as $a$.
Indeed, let $\alpha$ be the semilattice component of
$a$ and $\beta$ that of such a witness $e$.
Choose all testing variables from the $\alpha$-component.
The indices of the two evaluated word values are then
$\alpha$ and $\alpha\wedge\beta$, so equality gives
$\alpha\leq\beta$.  Choose them instead from the
$\beta$-component; equality gives
$\beta\leq\alpha$. Therefore $\alpha=\beta$.
These comparisons apply separately to every witness.

Consider first $\omega_1=x^{-1}y^{-1}xy$.
Let $e$ be a witness for \emph{left insertion in the first}
variable and $f$ a witness for \emph{right insertion in the
second} variable.  Their marginal tests include
\begin{equation}\label{eq:B-exist-A-tests}
 \omega_1(at,y)=\omega_1(et,y),\qquad
 \omega_1(y,ta)=\omega_1(y,tf)
 \quad(t,y\in S).
\end{equation}
Choose $t,y$ in the completely simple component
of $a$.  The row of a value of $\omega_1(x,y)$
comes from $x$, and its column comes from $y$.
The first equality in \eqref{eq:B-exist-A-tests}
therefore forces $e\mathrel{\mathcal R}a$;
the second forces $f\mathrel{\mathcal L}a$.
As $e,f,p$ are idempotents in that component, we have
\[
 e\mathrel{\mathcal R}p\ \Longrightarrow\ ep=p,
 \qquad
 f\mathrel{\mathcal L}p\ \Longrightarrow\ pf=p.
\]
Now set $t=p$ in both equations of
\eqref{eq:B-exist-A-tests}.  Since $ap=pa=a$,
we obtain
\[
 \omega_1(a,y)=\omega_1(p,y),
 \qquad \omega_1(y,a)=\omega_1(y,p)
 \quad(y\in S).
\]
Thus $a\mathrel{\eta_1}p$, and
Theorem~\ref{thm:B-pointwise} implies
$a\mathrel{\mathcal C}p$, or $a\in Z(pSp)$.

Now consider $\omega_2=(yx)^{-1}xy$.
Choose just one witness $e$ for \emph{left insertion in
the second} variable.  The corresponding condition is
\begin{equation}\label{eq:B-exist-B-test}
 \omega_2(y,at)=\omega_2(y,et)
 \quad(t,y\in S).
\end{equation}
As before, the witness belongs to the same component
as $a$. In that component, $\omega_2(x,y)$ has both
row and column of its second argument.  The row
comparison in \eqref{eq:B-exist-B-test} gives
$e\mathrel{\mathcal R}a$, hence $ep=p$.
Taking $t=p$ gives
\[
 \omega_2(y,a)=\omega_2(y,p)\qquad(y\in S).
\]
Restrict $y$ to the completely regular local monoid
$T=pSp$.  Here $a$ is a unit and $p$ is the identity.
The single-direction assertion of
Lemma~\ref{lem:B-pointwise-local} proves
$a\in Z(T)=Z(pSp)$.
Thus in either case $a\in K$; notably, this conclusion
was obtained without identifying the witnesses $e,f$
with $p$ at the outset.

Conversely, suppose $a\in K$ and let $p=a^0$.
Then $a\in Z(pSp)\cap H_p$, hence
$a\mathrel{\mathcal C}p$ by \eqref{eq:C}.
Theorem~\ref{thm:main} gives $a\sim_{\omega_i}p$,
so all left and right marginal identities are true
with the single idempotent witness $p$ for each position.
Therefore $a\in M_{\omega_i}(S)$ and also
$a\in\widetilde M_{\omega_i}(S)$.
Together with the canonical characterization in
Corollary~\ref{cor:common-margin}, this proves every
equality in \eqref{eq:B-existential-equals}.
\end{proof}

\noindent {\bf Acknowledgment:}   The authors acknowledge the help of AI tools  in preparing this manuscript.  The paper is supported partially  by the Nature Science Foundations of  China (12271442, 11661082).


\begin{thebibliography}{99}

\bibitem{Hall1940}
P.~Hall,
Verbal and marginal subgroups,
\emph{J. Reine Angew. Math.}
\textbf{182} (1940), 156--157.

\bibitem{Clifford1961}
A.~H.~Clifford and G.~B.~Preston,
\emph{The Algebraic Theory of Semigroups},
Vol.~I, Mathematical Surveys, No.~7,
American Mathematical Society,
Providence, RI, 1961.

\bibitem{Petrich1999}
M.~Petrich and N.~R.~Reilly,
\emph{Completely Regular Semigroups},
John Wiley \& Sons, New York, 1999.

\bibitem{Radovic2025}
J.~Radovi\'c and N.~Mudrinski,
Commutators in completely simple semigroups,
\emph{Semigroup Forum}
\textbf{110} (2025), 216--228.

\bibitem{Araujo2025}
G.~Ara\'ujo, J.~Ara\'ujo and M.~Kinyon,
Marginal subsemigroups and commutators
in inverse semigroups,
\emph{Semigroup Forum}
\textbf{111} (2025), 351--361.

\bibitem{Kinyon2026}
M.~Kinyon and D.~Stanovsk\'y,
Abelianness and centrality in inverse semigroups,
\emph{Proc. Edinburgh Math. Soc.}
\textbf{69} (2026), 160--176.

\end{thebibliography}
\end{document}